# THE "FUNDAMENTAL THEOREM" FOR THE ALGEBRAIC $K$-THEORY OF SPACES. III. THE NIL-TERM

JOHN R. KLEIN AND E. BRUCE WILLIAMS

ABSTRACT. In this paper we identify the "nil-terms" for Waldhausen's algebraic $K$-theory of spaces functor as the reduced $K$-theory of a category of equivariant spaces equipped with a homotopically nilpotent endomorphism.

## 1. INTRODUCTION

This is the third in a series of papers which concerns the decomposition

$$A^{fd}(X \times S^1) \simeq A^{fd}(X) \times \mathcal{B}A^{fd}(X) \times N_- A^{fd}(X) \times N_+ A^{fd}(X).$$

Here, $A^{fd}(X)$ is Waldhausen's algebraic $K$-theory of the space $X$ and $\mathcal{B}A^{fd}(X)$ is a certain non-connective delooping of it. The remaining factors on the right, called "nil-terms," are homotopy equivalent $[H_+],[H_+2]$. They have not been given a $K$-theoretic description thus far.

In this installment, we will identify the the nil-terms as a shifted copy of the reduced $K$-theory of a category whose objects are equivariant spaces equipped with a homotopically nilpotent endomorphism.

Let $X$ be a connected based space. Let $G.$ denote the Kan loop group of the total singular complex of $X$, and let $G$ denote the geometric realization of $G.$. Then the classifying space $BG$ has the weak homotopy type of $X$.

Define a category nil($X$) in which an *object* consists of a pair

$$(Y, f)$$

such that $Y$ is a based space with $G$-action and $f: Y \to Y$ is an equivariant map which is *homotopically nilpotent* under composition. Aditionally, we assume that $Y$ admits the structure of a based $G$-cell complex in which the action of $G$ is free away from the basepoint. A *morphism* $(Y, f) \to (Z, g)$ is a based $G$-map $e: Y \to Z$ such that $g \circ e = e \circ f$.

There is a full subcategory $\text{nil}_{fd}(X)$ of $\text{nil}(X)$ whose objects are those $Y$ which are *finitely dominated* in the sense that $Y$ is a retract up to homotopy of an object which is built up from a point by attaching a finite number of free $G$-cells. A morphism of $\text{nil}_{fd}(X)$ is a *weak equivalence* if and only if its underlying map of topological spaces is a weak homotopy equivalence. It is a *cofibration* if its underlying map of spaces is obtained up to isomorphism by attaching free $G$-cells.

With the above structure, it turns out that $\text{nil}_{fd}(X)$ is a category with cofibrations and weak equivalences. It therefore has a $K$-theory, which is denoted $K^{fd}(\text{nil}(X))$. The forgetful map $(Y, f) \mapsto Y$ gives rise to a map on $K$-theories

$$K^{fd}(\text{nil}(X)) \to A^{fd}(X).$$

Let $\widetilde{K}^{fd}(\text{nil}(X))$ denote its homotopy fiber.

We now can state our main result, which establishes the other half of the "fundamental theorem" for $A^{fd}(X)$:

**Main Theorem.** *There is a homotopy equivalence of functors*

$$\widetilde{K}^{fd}(\text{nil}(X)) \quad \simeq \quad \Omega N_+ A^{fd}(X).$$

*Remark.* The above result is used in the paper [GKM], where it is shown that the homotopy groups of $N_+ A^{fd}(X)$ are either trivial or infinitely generated. Another result of that paper determines $p$-complete homotopy type of $N_+ A^{fd}(*)$ in degrees $\leq 4p - 7$, for $p$ an odd prime.

*Acknowledgements.* A early draft of this paper was circulated in the mid 1990s. A revived interest in understanding nil phenomena [*], [F] was an inducement for us to finally publish. The authors wish to thank the SFB 343 at Bielefeld University, and especially, Friedhelm Waldhausen, for providing the atmosphere enabling us to accomplish this research.

## 2. Preliminaries

In what follows, we assume that the reader is familiar with the material of [H$_+$].

The spaces in this paper are to be given the compactly generated topology. Products are taken in the compactly generated sense. Let $M.$ be a simplicial monoid, and let $M = |M.|$ denote its geometric realization. If $Y$ and $Z$ are (based, left) $M$-spaces, we say that a based $M$-map $Y \to Z$ is *weak equivalence* if (and only if) it is a weak homotopy equivalence of underlying topological spaces. Let $\mathbb{T}(M)$ denote the category of based $M$-spaces and based $M$-maps.

Recall from [H$_+$] that $\mathbb{C}(M)$ denotes the full subcategory of $\mathbb{T}(M)$ whose objects are based $M$-spaces which are cofibrant, in the sense that that they are built up from a point by cell attachments, in the partial ordering defined by dimension, where the *cell* of dimension $n$ is given by
$$D^n \times M$$
with action defined by left translation.

An object of $\mathbb{C}(M)$ is *finite* if it is built up from a point by finitely many cell attachments (up to isomorphism). An object of $\mathbb{C}(M)$ is said to be *homotopy finite* if there exists a weak equivalence to a finite object. An object of $\mathbb{C}(M)$ is said to be *finitely dominated* if it is an equivariant retract of a homotopy finite object. Let $\mathbb{C}_{fd}(M)$ denote the full subcategory of $\mathbb{C}(M)$ whose objects are finitely dominated.

A *cofibration* of $\mathbb{T}(M)$ (or its subcategories $\mathbb{C}(M)$, $\mathbb{C}_{fd}(M)$) is a morphism $X \to Y$ such that $Y$ is obtained from $X$ by a sequence of cell attachments, where $n$-cells are attached over cells of dimension $\leq n-1$.

We let $h\mathbb{C}_{fd}(M)$ denote the subcategory of $\mathbb{C}_{fd}(M)$ defined by the weak equivalences. With respect to these conventions, $\mathbb{C}_{fd}(M)$ has the structure of a category of cofibrations and weak equivalences, and its $K$-theory is
$$A^{fd}(*; M) \quad := \quad \Omega |h\mathcal{S}.\mathbb{C}_{fd}(M)|,$$
where the right side is the loop space of the geometric realization of Waldhausen's $\mathcal{S}.$-construction of $\mathbb{C}_{fd}(M)$ ([W, p. 330]). If $M$ is the realization of a simplicial group, then $A^{fd}(*; M)$ is one of the definitions of $A^{fd}(BM)$ (cf. [W, p. 379], [H$_+$, 1.6]).

The category $\mathrm{nil}_{fd}(X)$ has *objects* specified by pairs $(Y, f)$ with $Y \in \mathbb{C}_{fd}(G)$ and a $G$-map $f: Y \to Y$ which is homotopically nilpotent (under composition) through morphisms of $\mathbb{C}_{fd}(G)$, i.e., there exists a non-negative integer $k$ such that the $k$-fold composite $f^{\circ k}$ is equivariantly null homotopic. A *morphism* $(Y, f) \to (Z, g)$ is a map $e: Y \to Z$ such that $g \circ e = e \circ f$. A *cofibration* of $\mathrm{nil}_{fd}(X)$ is a morphism whose underlying map of $G$-spaces $Y \to Z$ is a cofibration of $\mathbb{C}_{fd}(G)$. A *weak equivalence* is a morphism whose underlying map of spaces is a weak homotopy equivalence.

**Lemma 2.1.** *With respect to the above conventions, $\mathrm{nil}_{fd}(X)$ is a category with cofibrations and weak equivalences.*

*Proof.* The non-trivial thing to be verified is that the cobase change axiom holds. Given a diagram
$$(B, f_1) \leftarrow (A, f_0) \rightarrowtail (C, f_2)$$

we define the pushout to be $(B \cup_A C, f)$, where $f$ denotes $f_1 \cup_{f_0} f_2$. Choose a positive integer $k$ such that $f_i^{\circ k}$ is null homotopic, for $i = 0, 1, 2$. It will be sufficient to check that $f^{\circ k}$ is homotopically nilpotent. Let us rename $g_i = f_i^{\circ k}$ and $g = f^{\circ k}$. Then there is a commutative diagram up to equivariant homotopy

$$\begin{array}{ccccc} B \vee C & \longrightarrow & B \cup_A C & \xrightarrow{\delta} & \Sigma A \\ {\scriptstyle g_1 \vee g_2} \downarrow & & {\scriptstyle g} \downarrow & & \downarrow {\scriptstyle \Sigma g_0} \\ B \vee C & \longrightarrow & B \cup_A C & \xrightarrow[\delta]{} & \Sigma A \end{array},$$

where $\delta$ is the boundary map in the Barratt-Puppe seqence. Since $g_1$ and $g_2$ are null homotopic, it follows that $g$ may be expressed as $\gamma \circ \delta$ up to homotopy, for some map $\gamma \colon \Sigma A \to B \cup_A C$. It follows that $g^{\circ 2} \simeq \gamma \circ \delta \circ \gamma \circ \delta$ is equivariantly null homotopic, for $\delta \circ \gamma \circ \delta \simeq \delta \circ g$ is homotopic to $(\Sigma g_0) \circ \delta$. □

## 3. Another look at the projective line

Let $\mathbb{N}_-$ denote monoid of negative integers with generator $t^{-1}$ and $\mathbb{N}_+$ denote the monoid of positive integers with generator $t$. Let $G$ be the realization of a simplicial group $G_\cdot$.

Recall that the *mapping telescope* of an object $Y_+ \in \mathbb{C}_{fd}(G \times \mathbb{N}_+)$ is the object $Y_+(t^{-1}) \in \mathbb{C}_{fd}(G \times \mathbb{N}_+)$ defined by taking the categorical colimit of the sequence

$$\cdots \xrightarrow{t} Y_+ \xrightarrow{t} Y_+ \xrightarrow{t} \cdots.$$

Similarly, if $Y_- \in \mathbb{C}_{fd}(G \times \mathbb{N}_-)$ is an object, we have a mapping telescope $Y_+(t)$ given by the colimit of

$$\cdots \xrightarrow{t^{-1}} Y_+ \xrightarrow{t^{-1}} Y_+ \xrightarrow{t^{-1}} \cdots.$$

Define $\mathbb{D}_{fd}(G \times \mathbb{Z})$ to be the category whose *objects* are diagrams

$$Y_- \to Y \to Y_+$$

in which $Y_- \in \mathbb{C}_{fd}(G \times \mathbb{N}_-)$, $Y \in \mathbb{C}_{fd}(G \times \mathbb{Z})$ and $Y_+ \in \mathbb{C}_{fd}(G \times \mathbb{N}_+)$, and where the maps $Y_- \to Y$ and $Y_+ \to Y$ are required to be based and equivariant. Moreover, the induced morphisms

$$Y_-(t) \to Y(t) \cong Y \quad \text{and} \quad Y_+(t^{-1}) \to Y(t^{-1}) \cong Y$$

are required to be cofibrations. We take the liberty of specifying the object as a diagram or as a triple $(Y_-, Y, Y_+)$.

A *morphism* $(Y_-, Y, Y_+) \to (Z_-, Z, Z_+)$ of $\mathbb{D}_{fd}(G \times \mathbb{Z})$ is a morphism $Y_- \to Z_-$, a morphism $Y \to Z$ and a morphism $Y_+ \to Z_+$ so that the evident diagram commutes. A *cofibration* is a morphism in which the induced maps

$$Y \cup_{Y_-(t)} Z_-(t) \to Z \quad \text{and} \quad Y \cup_{Y_+(t^{-1})} Z_+(t^{-1}) \to Z$$

are cofibrations.

The *projective line* $\mathbb{P}_{fd}(G)$ of [H$_+$] is given by the full subcategory of $\mathbb{D}_{fd}(G \times \mathbb{Z})$ whose objects $(Y_-, Y, Y_+)$ satisfy an auxiliary condition, viz., that the induced maps $Y_-(t) \to Y$ and $Y_+(t^{-1}) \to Y$ are weak homotopy equivalences. A *cofibration* is a morphism which is a cofibration of $\mathbb{D}_{fd}(G \times \mathbb{Z})$. A *weak equivalence* is a morphism in which $Y_- \to Z_-$, $Y \to Z$ and $Y_+ \to Z_+$ are weak homotopy equivalences of spaces.

Let $\mathbb{D}_{fd}(G \times \mathbb{N}_-) \subset \mathbb{D}_{fd}(G \times \mathbb{Z})$ denote the full subcategory whose objects $(Y_-, Y, Y_+)$ satisfy the condition that $Y_-(t) \to Y$ is a weak equivalence. Similarly, define $\mathbb{D}_{fd}(G \times \mathbb{N}_+)$ to be the full subcategory whose objects $(Y_-, Y, Y_+)$ satisfy the condition that $Y_+(t^{-1}) \to Y$ is a weak equivalence.

A morphism $Y_-, Y, Y_+) \to (Z_-, Z, Z_+)$ of $\mathbb{D}_{fd}(G \times \mathbb{N}_+)$ is a *weak equivalence* if the map $Y_+ \to Z_+$ is a weak homotopy equivalence. It is a *cofibration* if it is so when considered in $\mathbb{D}_{fd}(G \times \mathbb{Z})$

Let $\mathbb{P}_{fd}^{h_{\mathbb{N}_+}}(G) \subset \mathbb{P}_{fd}(G)$ denote the full subcategory with objects $(Y_-, Y, Y_+)$ such that $Y_+$ acyclic.

**Proposition 3.1.** *There is a homotopy fiber sequence*

$$\Omega|h\mathcal{S}.\mathbb{P}_{fd}^{h_{\mathbb{N}_+}}(G)| \to \Omega|h\mathcal{S}.\mathbb{P}_{fd}(G)| \to \Omega|h\mathcal{S}.\mathbb{D}_{fd}(G \times \mathbb{N}_+)|\,.$$

*Proof.* Define a courser notion of weak equivalence on the projective line by specifying a morphism $(Y_-, Y, Y_+) \to (Z_-, Z, Z_+)$ to be an $h_{\mathbb{N}_+}$-*equivalence* if (and only if) the map $Y_+ \to Z_+$ is a weak equivalence. Application of the *fibration theorem* [W, 1.6.5] shows that the sequence

$$\Omega|h\mathcal{S}.\mathbb{P}_{fd}^{h_{\mathbb{N}_+}}(G)| \to \Omega|h\mathcal{S}.\mathbb{P}_{fd}(G)| \to \Omega|h_{\mathbb{N}_+}\mathcal{S}.\mathbb{P}_{fd}(G)|$$

is a fibration up to homotopy.

Let $\mathbb{P}_{fd}(G) \to \mathbb{D}_{fd}(G \times \mathbb{N}_+)$ denote the inclusion functor. By [H$_+$, §4] we have that the induced map

$$|h_{\mathbb{N}_+}\mathcal{S}.\mathbb{P}_{fd}(G)| \to |h\mathcal{S}.\mathbb{D}_{fd}(G \times \mathbb{N}_+)|$$

induces an isomorphism on homotopy groups in degrees $> 1$. Hence, the homotopy fiber of the induced map of loop spaces

$$\Omega|h_{\mathbb{N}_+}\mathcal{S}.\mathbb{P}_{fd}(G)| \to \Omega|h\mathcal{S}.\mathbb{D}_{fd}(G \times \mathbb{N}_+)|$$

is homotopically trivial.

It follows that the homotopy fiber of the map
$$\Omega|h\mathcal{S}.\mathbb{P}_{fd}(G)| \to \Omega|h\mathcal{S}.\mathbb{D}_{fd}(G \times \mathbb{N}_+)|$$
is identified with the homotopy fiber of the map
$$\Omega|h\mathcal{S}.\mathbb{P}_{fd}(G)| \to \Omega|h_{\mathbb{N}_+}\mathcal{S}.\mathbb{P}_{fd}(G)|.$$
The result follows. □

## 4. The 'characteristic sequence'

Let $(Y, f) \in \text{nil}_{fd}(X)$ be an object, and let $Y \otimes \mathbb{N}_- \in \mathbb{C}_{fd}(G)$ be the object given by
$$(Y \times \mathbb{N}_-)/(* \times \mathbb{N}_-).$$
Then $f$ induces a self-map of $Y \otimes \mathbb{N}_-$ which is given by $(y, r) \mapsto (f(y), r)$. We will denote this self-map also by $f$.

Let $Y_f$ be the *homotopy coequalizer* of the pair of maps
$$Y \otimes \mathbb{N}_- \underset{t^{-1}}{\overset{f}{\rightrightarrows}} Y \otimes \mathbb{N}_-,$$
where $t^{-1}$ denotes the map $(y, r) \mapsto (y, r-1)$. (Recall that the homotopy coequalizer of a pair of morphisms $\alpha, \beta \colon U \to V$ is defined to be the quotient of the disjoint union $V \amalg (U \times [0, 1])$ which is given by identifying $(u, 0)$ with $\alpha(u)$, $(u, 1)$ with $\beta(u)$ and $* \times [0, 1]$ with the basepoint of $V$.)

If we give $Y$ the structure of a based $(G \times \mathbb{N}_-)$-space by letting $\mathbb{N}_-$ act by means of $f$, then we also have a $(G \times \mathbb{N}_-)$-equivariant map
$$\pi_f \colon Y \otimes \mathbb{N}_- \to Y$$
which is given by $(y, r) \mapsto f^{-r}(y)$. Then $\pi_f$ coequalizes $f$ and $t^{-1}$, so by the universal property of the homotopy coequalizer, there is an induced map
$$Y_f \to Y,$$
which is $(G \times \mathbb{N}_-)$-equivariant.

**Lemma 4.1.** *The map $Y_f \to Y$ induces an isomorphism in reduced singular homology.*

*Proof.* Let $p \colon S^1 \to S^1 \vee S^1$ be the pinch map, and let $\rho \colon S^1 \to S^1$ be the reflection map. Then the composite
$$S^1 \xrightarrow{p} S^1 \vee S^1 \xrightarrow{\text{id} \vee \rho} S^1 \vee S^1$$
will be denoted $(1, -1)$.

The homotopy coequalizer induces a homotopy cofiber sequence
$$\Sigma(Y \otimes \mathbb{N}_-) \xrightarrow{t^{-1}-f} \Sigma(Y \otimes \mathbb{N}_-) \to \Sigma Y_f$$

where the first map is defined to be the composite

$$\Sigma(Y \otimes \mathbb{N}_-) \xrightarrow{(1,-1)\wedge \mathrm{id}} \Sigma(Y \otimes \mathbb{N}_-) \vee \Sigma(Y \otimes \mathbb{N}_-) \xrightarrow{t^{-1} \vee f} \Sigma(Y \otimes \mathbb{N}_-).$$

Taking reduced singular chains, we get an induced homotopy cofiber sequence of chain complexes

$$(1) \qquad C_*(Y) \otimes \mathbb{Z}[t^{-1}] \xrightarrow{t_*^{-1} - f_*} C_*(Y) \otimes \mathbb{Z}[t^{-1}] \longrightarrow C_*(Y_f).$$

Now, for any $\mathbb{Z}$-module $M$ equipped with self-map $f\colon M \to M$, we have an exact sequence of $\mathbb{Z}[t^{-1}]$-modules

$$(2) \qquad 0 \longrightarrow M \otimes \mathbb{Z}[t^{-1}] \xrightarrow{t^{-1} - f} M \otimes \mathbb{Z}[t^{-1}] \longrightarrow M_f \longrightarrow 0$$

in which $M_f$ denotes $M$ considered as a $\mathbb{Z}[t^{-1}]$-module where $t^{-1}$ acts via $f$ (see [B, p. 630]). This implies that the sequence (1) becomes exact when $C_*(Y_f)$ is replaced by $C_*(Y)$ by means of the chain map $C_*(Y_f) \to C_*(Y)$ which is induced by the map $Y_f \to Y$. Consequently, the five lemma implies that the chain map $C_*(Y_f) \to C_*(Y)$ is a quasi-isomorphism. □

*Remark* 4.2. The sequence (1) is a chain complex version of the so-called, "characteristic sequence" (2) of the module $M$. Consequently, it is not inappropriate to think of the homotopy coequalizer diagram

$$Y \otimes \mathbb{N}_- \underset{t^{-1}}{\overset{f}{\rightrightarrows}} Y \otimes \mathbb{N}_- \longrightarrow Y_f$$

as a kind of non-linear version of the characteristic sequence (of the object $Y$).

**Preliminary identification of $K(\mathrm{nil}_{fd}(X))$.** Define an exact functor

$$\mathrm{nil}_{fd}(X) \xrightarrow{\Phi} \mathbb{P}_{fd}^{h_{\mathbb{N}_+}}(G)$$

by

$$(Y, f) \mapsto (Y_f, Y_f(t), *),$$

where $Y_f$ is defined above.

In the other direction, define an exact functor

$$\mathbb{P}_{fd}^{h_{\mathbb{N}_+}}(G) \xrightarrow{\Psi} \mathrm{nil}_{fd}(X)$$

by

$$(Y_-, Y, Y_+) \mapsto (Y_-, t^{-1}).$$

To see that $\Psi$ is well-defined, let $(Y_-, Y, Y_+)$ be an object of $\mathbb{P}_{fd}^{h_{\mathbb{N}_+}}(G)$. Then $Y_+$ and $Y$ are acyclic. Hence $Y_-$ has an acyclic mapping telescope.

This implies that there exists a $k \in \mathbb{N}_-$ such that $t^k\colon Y_- \to Y_-$ is $G$-equivariantly null homotopic. Let $Z$ denote the quotient
$$Y_-/t^k(Y_-)$$
considered as an object of $\mathbb{C}(G)$. Then $Z$ is finitely dominated. This is a consequence of a cell-by-cell induction when $Y_-$ is a finite object of $\mathbb{C}(G \times \mathbb{N}_-)$. It therefore also true when $Y_-$ is is finitely dominated, since every finitely dominated obect of $\mathbb{C}(G \times \mathbb{N}_-)$ is a retract of a finite object up to homotopy, and the operation $Y_+ \mapsto Y_+/t^k(Y_-)$ is functorial. Since $t^k$ is $G$-equivariantly null homotopic, the identity map $Y_- \to Y_-$ factors through $Z$ up to homotopy. It follows that $Y_-$ is also a finitely dominated when considered as object of $\mathbb{C}(G)$. This shows that $(Y_-, t^{-1})$ is an object of $\text{nil}_{fd}(X)$.

**Lemma 4.3.** *The functors $\Psi$ and $\Phi$ induce mutually inverse homotopy equivalences on $K$-theory.*

*Proof.* The composite $\Psi \circ \Phi$ is given by
$$(Y, f) \quad \mapsto \quad (Y_f, t^{-1})$$
and 4.1 implies that there is a morphism $(Y_f, t^{-1}) \to (Y, f)$ which is a weak equivalence after taking a suitable number of suspensions. Since suspension induces a homotopy equivalence on the level of $K$-theory [W, 1.6.2], it follows that $\Psi \circ \Phi$ induces a homotopy equivalence.

The composite $\Phi \circ \Psi$ is given by
$$(Y_-, Y, Y_+) \quad \mapsto \quad (Y_-, Y_-(t), *)\,.$$
This admits an evident equivalence to the identity functor. Consequently $\Phi \circ \Psi$ induces a map which is homotopic to the identity on the level of $K$-theory. □

## 5. Proof of the main theorem

By 4.3, we have a homotopy equivalence,
$$\Omega|h\mathcal{S}.\text{nil}_{fd}(X)| \quad \simeq \quad \Omega|h\mathcal{S}.\mathbb{P}_{fd}^{h\mathbb{N}_+}(G)|\,.$$
Plugging this into 3.1, we obtain a homotopy fiber sequence
$$\Omega|h\mathcal{S}.\text{nil}_{fd}(X)| \to \Omega|h\mathcal{S}.\mathbb{P}_{fd}(G)| \to \Omega|h\mathcal{S}.\mathbb{D}_{fd}(G \times \mathbb{N}_+)|\,.$$
Let $\epsilon\colon \Omega|h\mathcal{S}.\mathbb{D}_{fd}(G \times \mathbb{N}_+)| \to \Omega|h\mathcal{S}.\mathbb{C}_{fd}(G)|$ denote the *augmentation* map of [H$_+$, 7.1], which is induced by
$$(Y_-, Y, Y_+) \quad \mapsto \quad Y/\mathbb{Z}\,,$$
where $Y/\mathbb{Z}$ denotes the orbit space under the $\mathbb{Z}$-action. Recall that the *nil-term* $N_+ A^{fd}(X)$ was defined to be the homotopy fiber of $\epsilon$.

Similarly, $\epsilon$ restricts to a map on $\Omega|h\mathcal{S}.\mathbb{P}_{fd}(G)|$. Denote the homotopy fiber of this restriction by $\Omega|h\mathcal{S}.\mathbb{P}_{fd}(G)|^\epsilon$. Consequently, we have an induced homotopy fiber sequence

$$\Omega|h\mathcal{S}.\mathrm{nil}_{fd}(X)| \to \Omega|h\mathcal{S}.\mathbb{P}_{fd}(G)|^\epsilon \to N_+ A^{fd}(X).$$

In was shown in [H$_+$, 7.6] that the second of these maps

$$\Omega|h\mathcal{S}.\mathbb{P}_{fd}(G)|^\epsilon \to N_+ A^{fd}(X)$$

is null homotopic. Moreover, it was shown in [H$_+$, 7.5] that there is a homotopy equivalence

$$\Omega|h\mathcal{S}.\mathbb{P}_{fd}(G)|^\epsilon \quad \simeq \quad \Omega|h\mathcal{S}.\mathbb{C}_{fd}(G)|$$

induced by the *global sections* functor $\Gamma\colon \mathbb{P}_{fd}(G) \to \mathbb{C}_{fd}(G)$ defined by

$$(Y_-, Y, Y_+) \quad \mapsto \quad CY_- \cup Y \cup CY_+,$$

where $CY_-$ denotes the cone on $Y_-$.

Assembling this information, we have a homotopy fiber sequence

$$(3) \qquad \Omega|h\mathcal{S}.\mathrm{nil}_{fd}(X)| \xrightarrow{\alpha} \Omega|h\mathcal{S}.\mathbb{C}_{fd}(G)| \xrightarrow{\beta} N_+ A^{fd}(X)$$

where $\alpha$ is induced by the functor $(Z, f) \mapsto \Sigma Z$ and $\beta$ is null homotopic. Since the suspension functor $\Sigma\colon \mathbb{C}_{fd}(G) \to \mathbb{C}_{fd}(G)$ induces a homotopy equivalence (by [W, 1.6.2]), we see that the homotopy fiber of $\alpha$ is homotopy equivalent to the homotopy fiber of the map $\alpha'$ which is induced by the forgetful map $(Z, f) \mapsto Z$.

On the one hand, the homotopy fiber of $\alpha'$ is $\widetilde{K}^{fd}(\mathrm{nil}(X))$, by definition. On the other hand, the homotopy fiber sequence (3) implies that the homotopy fiber of $\alpha$ is homotopy equivalent to $\Omega N_+ A^{fd}(X)$. We conclude that there is a homotopy equivalence

$$\widetilde{K}^{fd}(\mathrm{nil}(X)) \quad \simeq \quad \Omega N_+ A^{fd}(X).$$

This completes the proof of the theorem.

Wayne State University, Detroit, MI 48202
*E-mail address*: `klein@math.wayne.edu`

University of Notre Dame, Notre Dame, IN 46556
*E-mail address*: `williams.4@nd.edu`